\documentclass[a4paper,11pt]{article}
\usepackage{amssymb}
\usepackage{latexsym}
\newtheorem{theorem}{Theorem}[section]
\newtheorem{lemma}[theorem]{Lemma}

\newtheorem{definition}[theorem]{Definition}

\newtheorem{corollary}[theorem]{Corollary}
\newtheorem{remark}[theorem]{Remark}

\newtheorem{thm}[theorem]{Theorem}

\newcommand{\comment}[1]{}  

\newcommand{\res}{\restriction}

\newcommand{\CZF}{{\mathbf{CZF}}}

\newcommand{\RDC}{{\mathbf{RDC}}}

\newcommand{\prf}{{\bf Proof: }}

\newcommand{\und}{\,\wedge\,}

\newcommand{\oder}{\,\vee\,}

\newcommand{\beqs}{\begin{eqnarray*}}
\newcommand{\eeqs}{\end{eqnarray*}}

\newcommand{\BI}{{\mathbf{BI}}}

\newcommand{\beq}{\begin{eqnarray}}
\newcommand{\eeq}{\end{eqnarray}}

\newcommand{\paar}[1]{\langle #1\rangle}
\newcommand{\UU}{{\mathbf U}}
\newcommand{\VV}{{\mathbf V}}

\newcommand{\MM}{{\mathbf M}}

\newcommand{\CST}{{\mathbf{CST}}}

\newcommand{\natt}{\mathrm{nat}}

\newcommand{\bN}{{\mathbb N}}

\newcommand{\N}{{\mathbb N}}

\newcommand{\PAs}[1]
{#1^+}

\newcommand{\bes}{\begin{eqnarray*}}
\newcommand{\ees}{\end{eqnarray*}}

\newcommand{\rel}{\Vdash}

\newcommand{\LPO}{{\mathbf{LPO}}}
 
 \newcommand{\LLPO}{{\mathbf{LLPO}}}
\newcommand{\ACA}{{\mathbf{ACA}}}

\newcommand{\SAC}{\Sigma^1_1\mbox{-}{\mathbf{AC}}}

\newcommand{\EE}{\,\dot{\in}\,} 
\newcommand{\EEE}{\mathrm{EL}}
\newcommand{\NEEE}{\mathrm{NEL}}
\newcommand{\NEE}{\,\dot{\notin}\,} 

\newcommand{\Prog}{{\mathrm{Prog}}}

\newcommand{\WF}{{\mathrm{WF}}}
\newcommand{\TI}{{\mathrm{TI}}}
\newcommand{\II}{{\mathrm{I}}}
\newcommand{\Comp}{{\mathrm{Comp^E}}}
\newcommand{\Acomp}{{\mathrm{A_{\!_E}}}}
\newcommand{\MMo}{{\mathfrak{M}}}
\newcommand{\pl}{{\mathrm{pl}}}
\newcommand{\komp}[3]{\{#1\}^E(#2)\simeq #3}
\newcommand{\kompp}[2]{\{#1\}^E(#2)}
\newcommand{\gl}{\dot{=}}

\def\provx#1#2#3#4{
\setbox1=\hbox{\kern1.5pt$\scriptstyle#3$}
\def\zeichen{#2}
\ifx\zeichen\empty\setbox0=\hbox to .75em{}\else\setbox0=\hbox
{\kern1.5pt$\scriptstyle#2$}\fi
\dimen1=\dp0 \ifdim \dimen1=0pt
\advance \dimen1 by 1.5ex \else \advance \dimen1 by 1.2ex
\fi\dimen3=2ex\dimen4=.5ex\ifdim \wd0<\wd1 \dimen2=\wd1 \else \dimen2=\wd0
\fi\hbox{$#1\hskip 5pt minus5pt\vrule height\dimen3
depth\dimen4\raise\dimen1\copy0\hskip-1\wd0 \lower\ht1
\copy1\hskip-1\wd1\vrule width\dimen2 height.7ex depth-.6ex\hskip3pt
minus1.5pt#4\hskip2pt plus2pt minus2pt$}}

\def\prov#1#2#3{
\setbox1=\hbox{\kern1.5pt$\scriptstyle#2$}
\def\zeichen{#1}
\ifx\zeichen\empty\setbox0=\hbox to .75em{}\else\setbox0=\hbox
{\kern1.5pt$\scriptstyle#1$}\fi
\dimen1=\dp0
\ifdim \dimen1=0pt
\advance \dimen1 by 1.5ex \else \advance \dimen1 by 1.2ex
\fi\dimen3=2ex\dimen4=.5ex\ifdim \wd0<\wd1 \dimen2=\wd1 \else \dimen2=\wd0
\fi\hbox{\hskip0pt plus 4pt
$\vrule height\dimen3
depth\dimen4\raise\dimen1\copy0\hskip-1\wd0
\lower\ht1\copy1\hskip-1\wd1\vrule width\dimen2 height.7ex depth-.6ex
\hskip3pt minus1.5pt#3\hskip2pt plus2pt minus2pt$}}

\def\prv#1#2{
\setbox1=\hbox{\kern1.5pt$\scriptstyle#2$}
\ifx\zeichen\empty\setbox0=\hbox to .75em{}\else\setbox0=\hbox
{\kern1.5pt$\scriptstyle#1$}\fi
\dimen1=\dp0 \ifdim \dimen1=0pt
\advance \dimen1 by 1.5ex \else \advance \dimen1 by 1.2ex
\fi\dimen3=2ex\dimen4=.5ex\ifdim \wd0<\wd1 \dimen2=\wd1 \else \dimen2=\wd0
\fi\hbox{\hskip.5em$\vrule height\dimen3
depth\dimen4\raise\dimen1\copy0\hskip-1\wd0
\lower\ht1\copy1\hskip-1\wd1\vrule width\dimen2 height.7ex depth-.6ex
\hskip3pt minus1.5pt$}}

\mathchardef\str='1066
\def\negprov#1#2#3{
\setbox1=\hbox{\kern1.5pt$\scriptstyle#2$}
\setbox4=\hbox{$\str$}
\def\zeichen{#1}
\ifx\zeichen\empty\setbox0=\hbox to 1em{}\else\setbox0=\hbox
{\kern1.5pt$\scriptstyle#1$}\fi
\dimen1=\dp0
\ifdim \dimen1=0pt
\advance \dimen1 by 1.5ex \else \advance \dimen1 by 1.2ex
\fi\dimen3=2ex\dimen4=.5ex\ifdim \wd0<\wd1 \dimen2=\wd1 \else \dimen2=\wd0
\fi\hbox{\hskip.5em$\kern-1.9pt\raise1pt\copy4\kern-\wd4\kern1.9pt\vrule height\dimen3
depth\dimen4\raise\dimen1\copy0\hskip-1\wd0
\lower\ht1\copy1\hskip-1\wd1\vrule width\dimen2 height.7ex depth-.6ex
\hskip3pt minus1.5pt#3\hskip2pt plus2pt minus2pt$}}

\def\goed#1{\setbox5=\hbox{$#1$}\dimen1=.25em \dimen2=\dimen1 \advance \dimen2
by -1pt\hbox{\raise.65\ht5 \hbox{\vrule height.5\ht5 depth0pt width.4pt\vrule
height.5\ht5 width\dimen1 depth-.48\ht5}\kern-\dimen2\copy5\kern-\dimen2
\raise.65\ht5 \hbox{\vrule height .5\ht5 width\dimen1 depth-.48\ht5\vrule
height.5\ht5 depth 0pt width.4pt}\hskip4pt plus2pt minus2pt}}

\def\mod#1#2{
\def\zeichen{#1}
\hbox{\hskip 2pt plus3pt minus 2pt\vrule width.5pt height2ex depth.5ex
\vbox{\ifx\zeichen\empty\hbox to .75em{}\else
\hbox{\kern1.5pt $\scriptstyle#1$}\fi
\kern2pt
\hrule
\kern1.7pt
\hrule\kern1.7pt}
\hskip3pt minus 2pt$#2$}\hskip2pt
plus3pt minus2pt}

\def\notmod#1#2{\hbox{\hskip 2pt plus 3pt minus 3pt\vrule width.5pt
height2ex depth.5ex
\vbox{\hbox{\kern1.5pt $\scriptstyle#1$}\kern3pt
\setbox0=\hbox{\kern2pt$\scriptstyle/$}
\hrule
\kern-1.7pt
\copy0
\kern-\ht0
\kern 1.7pt
\hrule\kern1.7pt}n
\hskip3pt minus 2pt$#2$}\hskip2pt
plus3pt minus2pt}
\def\sq{\hbox{\rlap{$\sqcap$}$\sqcup$}}
\def\qed{\ifmmode\sq\else{\unskip\nobreak\hfil\penalty50\hskip1em\null
\nobreak\hfil\sq\parfillskip=0pt\finalhyphendemerits=0\endgraf}\fi\medskip}

\def\lleq{\hbox{\hskip3pt minus3pt\kern1pt\lower4pt
\vbox{\hbox{$\scriptstyle\ll$}
\kern-7pt\hbox{\kern1pt$\scriptstyle=$}}\hskip3pt minus 3pt}}

\mathchardef\res='1152
\mathchardef\qin='1062
\mathchardef\qprec='1036
\mathchardef\qless='474
\mathchardef\dpkt='72

\begin{document}

\title{Constructive Zermelo-Fraenkel set theory and the limited principle of omniscience}

\author{Michael Rathjen\\
Department of Pure Mathematics\\
University of Leeds, Leeds LS2 9JT, England\\  E-mail:~{\sf 
rathjen@maths.leeds.ac.uk} }
\maketitle
\begin{abstract}
In recent years the question of whether adding the limited principle of omniscience, $\LPO$,
to constructive Zermelo-Fraenkel set theory, $\CZF$, increases its strength has arisen several times.
As the addition of excluded middle for atomic formulae to $\CZF$ results in a rather strong theory,
i.e. much stronger than classical Zermelo set theory, it is not obvious that its augmentation by $\LPO$ would be proof-theoretically
benign. The purpose of this paper is to show that $\CZF+\RDC+\LPO$ has indeed the same strength as $\CZF$, where $\RDC$ stands for relativized dependent choice. In particular, these theories prove the same $\Pi^0_2$ theorems of arithmetic.
\\[1ex]
Keywords: Constructive set theory, limited principle of omniscience,
bar induction,
proof-theoretic strength \\ MSC2000: 03F50;  03F25;
03E55;  03B15; 03C70

\end{abstract}

\section{Introduction}
Constructive Set Theory was introduced by
 John Myhill in a seminal paper \cite{myhill}, where a
specific axiom system $\CST$ was introduced.
 Through developing constructive
 set theory he wanted to isolate the principles underlying Bishop's conception
of what sets and functions are, and he  wanted  ``these
principles to be such as to
 make the process of formalization completely trivial,
 as it is in the classical case"
(\cite{myhill}, p. 347).
 Myhill's $\CST$ was subsequently modified by Aczel
 and  the resulting
  theory was called {\em Zermelo-Fraenkel set theory},
  $\CZF$. A hallmark of this theory is that it possesses a type-theoretic
interpretation (cf. \cite{aczel82,mar}).
Specifically, $\CZF$ has a scheme called Subset Collection Axiom
(which is a generalization of Myhill's Exponentiation Axiom) whose
formalization was directly inspired by the type-theoretic
interpretation.

Certain basic principles of classical mathematics are taboo for
the constructive mathematician. Bishop called them {\em principles
of omniscience}.
The limited principle of omniscience, $\LPO$, is an instance of the law of excluded middle which usually
serves as a line of demarcation, separating ``constructive" from ``non-constructive" theories.
Over the last few years the question of whether adding  $\LPO$
to constructive Zermelo-Fraenkel set theory increases its strength has arisen several times.
As the addition of excluded middle for atomic formulae to $\CZF$ results in a rather strong theory,
i.e. much stronger than classical Zermelo set theory, it is not obvious that its augmentation by $\LPO$ would be proof-theoretically
benign. The purpose of this paper is to show that $\CZF+\RDC+\LPO$ has indeed the same strength as $\CZF$, where $\RDC$ stands for relativized dependent choice. In particular, these theories prove the same $\Pi^0_2$ theorems of arithmetic. The main tool will be a realizability model for $\CZF+\RDC+\LPO$ that is based on recursion in a type-2 object.
This realizability interpretation is shown to be formalizable in the theory of bar induction, $\BI$, which is known to have the same strength as $\CZF$.

To begin with we recall some principles of omniscience.
Let $2^{\bN}$ be Cantor space, i.e the set of all functions from the naturals into $\{0,1\}$.
 \begin{definition}\label{omnis} {\em
 { Limited Principle of Omniscience} ($\LPO$):
$$\forall f\in 2^{\bN}\,[\exists n\,f(n)=1\;\;\vee\;\;\forall n\,f(n)=0].$$
 { Lesser Limited Principle of Omniscience} ($\LLPO$):
\begin{eqnarray*}\forall f\in 2^{\bN}\,
\bigl(\forall n,m[f(n)=f(m)=1\to n=m]
\,\to\,[\forall n\,f(2n)=0\;\;\vee\;\;
 \forall n\,f(2n+1)=0]\bigr).\end{eqnarray*}
 }\end{definition}
 $\LPO$ is incompatible with both Brouwerian mathematics and
 Russian constructivism. With $\LLPO$ the story is more complicated
 as it is by and large compatible with Russian constructivism, namely with the form of Church Thesis saying that every function from
 naturals to naturals is computable (recursive) even on the basis of full intuitionistic Zermelo-Fraenkel set theory (see \cite{RMC}).

\section{The theory $\BI$}
In the presentation of subsystems of second order arithmetic we follow \cite{SOSA}. By $\mathcal{L}_2$ we denote the
language of these theories.
$\mathbf{ACA}_0$ denotes the theory of arithmetical comprehension.
\begin{definition}{\em
For a 2-place relation $\prec$ and an arbitrary formula $F(a)$ of $\mathcal{L}_2$ we define
\begin{enumerate}
\item[]
$\Prog(\prec,F):=\forall x[\forall y (y\prec x \rightarrow F(y))\rightarrow F(x)]$ (\emph{progressiveness})
\item[]$\TI(\prec,F):= \Prog(\prec,F)\rightarrow\forall x F(x)$ (\emph{transfinite induction})
\item[] $\WF(\prec):=\forall X\TI(\prec,X):=$ \ \newline
$\forall X(\forall x[\forall y (y\prec x \rightarrow y\in X))\rightarrow x\in X]\rightarrow \forall x[x\in X])$ (\emph{well-foundedness}).
\end{enumerate}
Let $\mathcal{F}$ be any collection of formul\ae \ of $\mathcal{L}_2$. For a 2-place relation $\prec$ we will write $\prec\in\mathcal{F}$, if $\prec$   is defined by a formula $Q(x,y)$ of $\mathcal{F}$ via $x\prec y:=Q(x,y)$.

The bar induction scheme is the collection of  all formul\ae \ of the form
$$\WF(\prec)\rightarrow\TI(\prec,F),$$
where $\prec$ is an arithmetical relation (set parameters allowed) and $F$ is an arbitrary formula of $\mathcal{L}_2$.

The theory $\ACA_0+\mbox{bar induction}$ will be denoted by $\BI$. In Simpson's book the acronym used for bar induction is $\Pi^1_{\infty}\mbox{-}\mathbf{TI}_0$
(cf. \cite[\S VII.2]{SOSA}).
}\end{definition}

\begin{thm} The following theories have the same proof-theoretic strength:
\begin{itemize}
\item[(i)] $\BI$
\item[(ii)] $\CZF$
\item[(iii)] The theory $\mathbf{ID}_1$ of (non-iterated) arithmetical inductive definitions.
\end{itemize}
\end{thm}

There is an interesting other way of characterizing $\BI$ which uses the notion of a countable coded $\omega$-model.

\begin{definition}{\em Let $T$ be a theory in the language of second order arithmetic, $\mbox{L}_2$. A {\em countable coded $\omega$-model of $T$} is a set $W\subseteq {\mathbb N}$, viewed as encoding the $\mbox{L}_2$-model
  $${\mathbb M}=({\mathbb N},{\mathcal S},+,\cdot,0,1,<)$$
  with ${\mathcal S}=\{(W)_n\mid n\in{\mathbb N}\}$ such that ${\mathbb M}\models T$
  (where $(W)_n=\{m\mid \paar{n,m}\in W\}$; $\paar{\,,}$ some coding function).

  This definition can be made in $\mathbf{RCA}_0$ (see \cite{SOSA}, Definition VII.2).

  We write $X\in W$ if $\exists n\;X=(W)_n$.
  }\end{definition}

  \begin{thm}\label{omega-ref} $\BI$ proves $\omega$-model reflection, i.e., for every formula $F(X_1,\ldots,X_n)$ with all free second order
  variables exhibited,
  $$\BI\vdash F(X_1,\ldots,X_n)\to \exists {\mathbb M}[\mbox{$\mathbb M$ countable coded $\omega$ model of $\ACA_0$ $\wedge$
  $\vec X\in {\mathbb M}$ $\wedge$ ${\mathbb M}\models F(\vec X\,)$.}]$$

\end{thm}
\prf \cite[Lemma VIII.5.2]{SOSA}. \qed
\begin{definition}{\em The scheme of $\SAC$ is the collection of all formulae
$$\forall x\,\exists X\,F(x,X)\to \exists Y\,\forall x\,F(x,(Y)_x)$$
with $F(x,X)$ of complexity $\Sigma^1_1$.

}\end{definition}

\begin{corollary}\label{SAC} $\BI$ proves that for every set $X$ there exists a countable coded $\omega$-model of $\ACA_0+\SAC$
containing $X$. In particular, $\BI$ proves $\SAC$ and $\Delta^1_1$-comprehension.
\end{corollary}
\prf Suppose $\forall x\,\exists X\,F(x,X,\vec U\,)$. Then there exists a countable coded $\omega$-model $\mathbb M=(\N,\mathcal S,+,\cdot,0,1,<)$
with $\vec U\in\mathbb M$ and $\mathbb M\models \forall x\,\exists X\,F(x,X,\vec U\,)$.
Let $\mathcal S=\{(W)_n\mid n\in\N\}$.
Define $f(n)=m$ if $\mathbb M\models F(n,(W)_m)$ and for all $k<m$ $\mathbb M\models \neg F(n,(W)_k)$.
Put $Y:=\{\langle n,x\rangle\mid x\in (W)_{f(n)}\}$. We then have $\mathbb M\models F(n,(Y)_n,\vec U\,)$ for all
$n$. Since $F$ is $\Sigma^1_1$ it follows that $ F(n,(Y)_n,\vec U\,)$ holds for all $n$.
This shows $\SAC$.

To show that for any set $Z$ there is an $\omega$-model of $\SAC$ containing $Z$ just
note that $\ACA_0+\SAC$ is finitely axiomatizable.

$\Delta^1_1$-comprehension is a consequence of $\SAC$. \qed

\begin{lemma}\label{All-Schluss} Let $A(X)$ be an arithmetic formula and $F(x)$ be an arbitrary formula of $\mathcal{L}_2$.
Let $A(F)$ be the formula that arises from $A(X)$ by replacing every subformula $t\in X$ by $F(t)$
(avoiding variable clashes, of course). Then we have
$$\BI\vdash \forall X\,A(X)\to A(F)\,.$$
\end{lemma}
\prf Arguing in $\BI$ suppose that $\neg A(F)$. Pick an $\omega$-model $\MM$ of $\ACA_0$ containing all parameters from $A$ and $F$ such that
$\MM\models \neg A(F)$. Letting $U=\{n\mid \MM\models F(n)\}$ we have
 $\neg A(U)$ because $A$ is an arithmetic formula and $\MM$ is absolute for such formulae on account of being an $\omega$-model.
 Thus we have shown $$\BI\vdash \neg A(F)\to \exists X\,\neg A(X)$$ from which the desired assertion follows. \qed

\section{Inductive definitions in $\BI$}
\begin{definition} Let $A(x,X)$ be an arithmetic formula in which the variable $X$ occurs positively. Henceforth we shall
convey this by writing $A(x,X^+)$.

Define \begin{eqnarray}\label{IA} \II_A(u) &:\Leftrightarrow & \forall X\,[\forall x\,(A(x,X)\to x\in X) \to u\in X]\,.\end{eqnarray}
We write $\II_A\subseteq F$ for $\forall v\,(\II_A(v)\to F(v))$, and $ F\subseteq \II_A$
for $\forall v\,(F(v) \to \II_A(v))$.
\end{definition}

\begin{lemma}\label{ind} The following are provable in $\BI$ for every $X$-positive arithmetic formula $A(x,X^+)$ and
arbitrary $\mathcal{L}_2$ formula $F(u)$.
\begin{itemize}
\item[(i)] $\forall u\,(A(u,\II_A)\to u\in\II_A)$.
\item[(ii)] $\forall x\,[A(x,F)\to F(x)] \to \II_A\subseteq F\,$ 
\item[(iii)]  $\forall u\,( u\in\II_A\to A(u,\II_A))$.
\end{itemize}
\end{lemma}
\prf (i): Assume $A(u,\II_A)$ and $\forall x\,(A(x,X)\to x\in X)$. The latter implies $\II_A\subseteq X$. Since $A(u,\II_A)$ holds,
and
owing to the positive occurrence
of $\II_A$ in the latter formula, we have $A(u,X)$. Since $X$ was arbitrary, we conclude that $\II_A(u)$.

(ii): Suppose $\II_A(u)$. Then $\forall X\,[\forall x\,(A(x,X)\to x\in X) \to u\in X]$, and hence, using Lemma \ref{All-Schluss},
$\forall x\,(A(x,F)\to F(x)) \to F(u)$. Thus, assuming $\forall x\,(A(x,F)\to F(x))$, we have $F(u)$.

(iii): Let $F(v):\Leftrightarrow A(v,\II_A)$. By (i) we have $F\subseteq \II_A$.
 Assuming $A(u,F)$ it therefore follows that $A(u,\II_A)$ since $F$ occurs positively in the former formula,
 and hence $F(u)$. Thus, in view of (ii), we get $\II_A\subseteq F$, confirming (iii).
\qed

\section{Recursion in a type-2 object}
Using the apparatus of inductive definitions, we would like to formalize in $\BI$ recursion in the type 2 object
$E:(\N\to \N)\to \N$ with $E(f)=n+1$ if $f(n)=0$ and $\forall i<n\,f(n)>0$ and $E(f)=0$ if $\forall n\,f(n)>0$.

In the formalization we basically follow \cite[VI.1.1]{hinman}.
We use some standard coding of tuples of natural numbers. The code of the empty tuple is $\langle\rangle:=1$, and
for any $k>0$ and tuple $(m_1,\ldots,m_k)$ let $\langle m_1,\ldots m_k\rangle:=p_1^{m_1+1}\cdot\ldots\cdot p_k^{m_k+1}$,
where $p_i$ denotes the $i$-th prime number.

\begin{definition}\label{E-recursion} Below $\mathrm{Sb}_0$ denotes the primitive recursive function from \cite[II.2.5]{hinman}
 required for what is traditionally called the S-m-n theorem.  Let $\Comp$ be the smallest class such that for all $k,n,p,r$, and $s$, all $i<k$
and ${\mathbf m}=m_1,\ldots,m_k$ in $\N$,
\begin{itemize}
\item[(0)] $\langle \langle 0,k,0,n\rangle,{\mathbf m},n)\in\Comp$;

$\langle \langle 0,k,1,i\rangle,{\mathbf m},m_i\rangle \in\Comp$;

$\langle \langle 0,k,2,i\rangle,{\mathbf m},m_i+1\rangle \in\Comp$;

$\langle \langle 0,k+3,4\rangle,p,q,r,s,{\mathbf m},p\rangle \in\Comp$ if $r=s$;

$\langle \langle 0,k+3,4\rangle,p,q,r,s,{\mathbf m},q\rangle \in\Comp$ if $r\ne s$;

$\langle \langle 0,k+2,5\rangle,p,q,{\mathbf m},{\mathrm{Sb}}_0(p,q\rangle \rangle \in\Comp$;

\item[(1)] for any $k',b,c_0,\ldots,c_{k'-1},q_0,\ldots,q_{k'-1}$ in $\N$, if for all $i<k'$
$\langle c_i,{\mathbf m},q_i\rangle \in\Comp$ and $\langle b,{\mathbf q},n\rangle \in\Comp$, then
$$\langle \langle 1,k,b,c_0,\ldots,c_{k'-1}\rangle, {\mathbf m},n\rangle \in\Comp\,;$$

\item[(2)] for any $b\in\N$, if $\langle b,{\mathbf m},n\rangle \in\Comp$, then
$$\langle \langle 2,k+1\rangle,b,{\mathbf m},n\rangle \in\Comp\,.$$

\item[(3.1)] for any $b\in\N$, if for all $p\in\N$ there exists $k_p\in\N$ with $k_p>0$
and $\langle b,p,{\mathbf m},k_p\rangle \in\Comp$, then
$$\langle \langle 3,k,b\rangle,{\mathbf m},0\rangle \in\Comp\,.$$
\item[(3.2)] for any $b,p\in\N$, if $\langle b,p,{\mathbf m},0\rangle \in\Comp$ and for all $i<p$ there exists $k_i\in\N$ with $k_i>0$
and $\langle b,i,{\mathbf m},k_i\rangle \in\Comp$, then
$$\langle \langle 3,k,b\rangle,{\mathbf m},p+1\rangle \in\Comp\,.$$
\end{itemize}

Clearly $\Comp$ is defined by a positive arithmetic inductive definition that we denote by $\Acomp$, i.e.,
$\Comp=\II_{\Acomp}$.
\end{definition}
\begin{lemma}\label{eindeut} For all $a,\mathbf m\in\N$ there is at most one $n\in\N$ such that $\langle a,\mathbf m,n\rangle \in\Comp$.
\end{lemma}
\prf Define a class $\mathfrak X$ by
 \begin{eqnarray*}\langle a,\mathbf m,n\rangle \in\mathfrak X &\mbox{ iff }& \langle a,\mathbf m,n\rangle \in\Comp\mbox{ and for all $k\in \N$, if $\langle a,\mathbf m,k\rangle \in\Comp$, then $n=k$.}\end{eqnarray*}

By Lemma  \ref{ind} (ii) we only have to show that $\mathfrak X$ is closed under the clauses defining $\Comp$. This is a straightforward affair, albeit a bit tedious. \qed

We shall put to use this notion of computability  for a realizability interpretation of $\CZF+\LPO$. This, however, will require
that the computability relation be a set  rather than a class such as $\Comp$.
To achieve this we shall invoke Theorem \ref{omega-ref}.

\begin{lemma} $\BI$ proves that there exists a countable coded $\omega$-model $\MMo$ of $\ACA$ such that the following hold.
\begin{itemize}
\item[(i)] $\MMo\models \forall x,\mathbf{y},z\,[\langle x,\mathbf{y},z\rangle \in\Comp \leftrightarrow \Acomp(\langle x,\mathbf{y},z\rangle ,\Comp)]$.
\item[(ii)] $\MMo\models  \forall x,\mathbf{y},z,z'[\langle x,\mathbf{y},z\rangle \in\Comp \und \langle x,\mathbf{y},z'\rangle \in\Comp \to z=z']$.
\end{itemize}
\end{lemma}
\prf This follows from Lemma \ref{ind} and Lemma \ref{eindeut} using Theorem \ref{omega-ref}. \qed

We will fix a model $\MMo$ as in the previous Lemma for the remainder of the paper and shall write
\begin{eqnarray*} \{a\}^E(\mathbf m)\simeq n & \Leftrightarrow & \MMo\models \langle a,\mathbf m,n\rangle \in\Comp\,.\end{eqnarray*}

Note that this notion of computability hinges on $\MMo$. More computations might converge in $\MMo$ than
outside of $\MMo$.

\section{Emulating a type structure in $\BI$}
We would like to define a type-theoretic interpretation of $\CZF+\RDC+\LPO$ in $\BI$.
 This will in a sense be similar to Aczel's interpretation of $\CZF$ in Martin-L\"of type theory (cf. \cite{aczel82}).
 To this end, we initiate a simultaneous  positive inductive definition of types $\UU$ and their elements as well as non-elements, and
 also a type $\VV$ of (of codes of) well-founded trees over $\UU$.
 The need for defining both elementhood and non-elementhood for types is necessitated by the requirement of positivity of the inductive definition.
 \\[1ex]
Below we use  the pairing function $\jmath(n,m) = (n+m)^2+n+1$
and its inverses
 $()_0,()_1$  satisfying
 $(\jmath(n,m))_0=n$ and $(\jmath(n,m))_1=m$. We will just write $(n,m)$ for $\jmath(n,m)$.

 \begin{definition}\label{Erec} Let $n_{\N}:=(0,n)$, $\natt:=(1,0)$, $\pl(n,m):=(2,(n,m))$, $\sigma(n,m):=(3,(n,m))$, $\pi(n,m):=(4,(n,m))$, and
 $\sup(n,m):=(5,(n,m))$.

 We inductively define classes $\UU,\EEE,\NEEE$ and $\VV$ 
 by the following clauses. Rather than $(n,m)\in \EEE$ and $(n,m)\in \NEEE$ 
  we write
 $n\EE m$ and $n\NEE m$, 
 respectively.
 \begin{enumerate}
 \item $n_{\N}\in \UU$; if $k<n$ then $k\EE n_{\N}$; if $k\geq n$ then $k\NEE n_{\N}$.
 \item $\natt \in\UU$ and  $n\EE \natt$ for all $n$.
 \item If $n,m\in\UU$, then $\pl(n,m)\in\UU$.
 \item Assume $\pl(n,m)\in \UU$.

 If $k\EE n$, then $(0,k)\EE\pl(n,m)$. If $k\EE m$, then $(1,k)\EE\pl(n,m)$.

 If $k\NEE n$, then $(0,k)\NEE\pl(n,m)$. If $k\NEE m$, then $(1,k)\NEE\pl(n,m)$. If $k$ is neither of the form $(0,x)$ nor
 $(1,x)$ for some $x$, then $k\NEE\pl(n,m)$.

 \item If $n\in \UU$ and $k\NEE n \oder \exists x\,( \komp ekx \und x\in \UU)$ holds for all $k$, then $\sigma(n,e)\in \UU$.
 \item  Assume $\sigma(n,e)\in \UU$.

 If $k\EE n$ and  $\exists x\,(\komp ekx\und u\EE x)$, then $(k,u)\EE \sigma(n,e)$.

 If $k\NEE n$ or $\exists x\,(\komp ekx\und u\NEE x)$, then $(k,u)\NEE  \sigma(n,e)$.

 If $x$ is not of the form $(u,v)$ for some $u,v$, then $x\NEE \sigma(n,e)$.

 \item If $n\in \UU$ and $k\NEE n \oder \exists x\,( \komp ekx \und x\in \UU)$ holds for all $k$, then $\pi(n,e)\in \UU$.
 \item  Assume $\pi(n,e)\in \UU$. 

 If $k\NEE n \oder \exists x,y\,(\komp ekx \und \komp dky\und y\EE x)$ holds for all $k$, then $d\EE\pi(n,e)$.

 If $\exists u \,(u\EE n \und \forall z\,\neg\komp duz)$, then $d\NEE\pi(n,e)$.

 If $\exists u \,\exists x\,(u\EE n \und \komp eux\und \exists z\,(\komp duz \und z\NEE x ))$, then $d\NEE\pi(n,e)$.

 \item If $n\in \UU$ and $k\NEE n\oder \exists x\,(\komp ekx\und x\in \VV)$ holds for all $k$, then
 $\sup(n,e)\in \VV$.

 \end{enumerate}
\end{definition}

\begin{remark}{\em Clearly, the predicates $\UU,\EE,\NEE$ and $\VV$ 
all appear positively in the above inductive definition.
Moreover, it falls under the scope of arithmeical inductive definitions and is therefore formalizable in our background
 theory $\BI$. Note also that for this it was important to move from the $\Pi^1_1$ computability notion of Definition
 \ref{E-recursion} to $E$-recursion in the $\omega$-model $\mathfrak M$.}\end{remark}

$\EE$ and $\NEE$ 
are complementary in the following sense.

\begin{lemma}\label{equality} For all $n\in \UU$, 
  \begin{eqnarray*}&&\forall x\,(x\EE n\leftrightarrow \neg\, x\NEE n)\,.
 \end{eqnarray*}
  \end{lemma}

  \prf This can be proved by the induction principle of Lemma \ref{ind}(ii). \qed

\begin{corollary}\label{set} For each $n\in \UU$,
  $\{x\mid x\EE n\}$ is a set.
  \end{corollary}

  \prf Note that $\EE$ and $\NEE$ are $\Pi^1_1$ as they are given by positive arithmetical inductive definitions.
  Since $\BI$ proves $\Delta^1_1$-comprehension by Corollary \ref{SAC}, it follows from Lemma \ref{equality}
  that  $\{x\mid x\EE n\}$ is a set.\qed

  \begin{definition} We shall use lower case Greek letters $\alpha,\beta,\gamma,\delta,\ldots$ to range over elements
  of $\VV$.

  Using the induction principle from Lemma \ref{ind}(ii), one readily shows that
  every $\alpha\in\VV$ is of the form $\sup(n,e)$ with $n\in \UU$ and $\forall x\EE n\;\kompp ex\in \VV$,
  where $\kompp ex\in \VV$ is an abbreviation for $\exists y\;(\komp exy\und y\in\VV)$.

  If $\alpha=\sup(n,e)$ we denote $n$ by $\bar{\alpha}$ and $e$ by $\tilde{\alpha}$.
  For $i\EE \bar{\alpha}$ we  shall denote by $\tilde{\alpha}i$  the unique $x$ such that $\komp {\tilde{\alpha}}ix$.
  \end{definition}

  If $\wp$ is an $r+1$-ary partial $E$-recursive function we denote by $\lambda x.\wp(x,\vec a\,)$ an index of the function
  $x\mapsto \wp(x,\vec a\,)$ (say provided by the S-m-n theorem or parameter theorem).

  \begin{lemma}\label{=}
  There is a 2-ary partial $E$-recursive function $\gl$ such that $\gl(\alpha,\beta)$ is defined for all $\alpha,\beta\in\VV$
 and (writing in infix notation $\alpha\gl\beta$ for $\gl(\alpha,\beta)$) the following equation holds
  \begin{eqnarray}\label{gl} (\alpha\gl\beta) & = & \sigma(\pi(\bar{\alpha},\lambda x.\sigma(\bar{\beta},\lambda y.(\tilde{\alpha} x\gl \tilde{\beta}y))),\lambda z. \pi(\bar{\beta},\lambda y.\sigma(\bar{\alpha},\lambda x.(\tilde{\alpha} x\gl\tilde{\beta}y))))\,.\end{eqnarray}
\end{lemma}
\prf Such a function can be defined by the recursion theorem for $E$-recursion. Totality on $\VV\times\VV$ follows from the
induction principle for $\VV$. \qed
  \section{Realizability}
  We will introduce a realizability semantics for sentences of set theory with
parameters from $\VV$. Bounded set quantifiers will be treated as
quantifiers in their own right, i.e., bounded and unbounded
quantifiers are treated as syntactically different kinds of
quantifiers. Let $\alpha,\beta\in \VV$ and $e,f\in \N$. We write
$e_{i,j}$ for $((e)_i)_j$.

  \begin{definition}[Kleene realizability over $\VV$]
Below variables $e,d$ range over natural numbers.
We define
\begin{eqnarray*}
e\rel\alpha=\beta &\mbox{iff}& e\EE (\alpha\gl \beta)\\
 e\rel\alpha\in \beta &\mbox{iff} & (e)_0\EE \bar{\beta}\;\wedge\;
(e)_1\rel\alpha=\tilde{\beta}(e)_0
 \\
e\rel\phi\wedge\psi &\mbox{iff}& (e)_0\rel\phi\;\wedge\;(e)_1\rel\psi\\
e\rel\phi\vee\psi &\mbox{iff}& \bigl[(e)_0={\mathbf
0}\,\wedge\,(e)_1\rel\phi\bigr]\;\vee
\;\bigl[(e)_0={\mathbf 1}\;\wedge\;(e)_1\rel\psi\bigr]\\
e\rel\neg\phi \phantom{AA} &\mbox{iff}& \forall d\;\neg d\rel\phi\\
e\rel\phi\rightarrow\psi &\mbox{iff}& \forall d\,
\bigl[d\rel\phi\;\rightarrow\;
\kompp ed\rel\psi \bigr]\\   
e\rel\forall x\in \alpha\; \phi(x)  &\mbox{iff}& \forall
i\EE\bar{\alpha}\; \kompp ei\rel\phi(\tilde{\alpha}i)
\\
e\rel\exists x\in \alpha\;\phi(x)  &\mbox{iff}&
(e)_0\EE\bar{\alpha}\;\wedge\; (e)_1\rel\phi(\tilde{\alpha}(e)_0)
 \\
e\rel\forall x \phi(x) \phantom{Ai} &\mbox{iff}& \forall \alpha\in \VV \;\kompp e{\alpha}\rel\phi(\alpha)\\
e\rel\exists x\phi(x)\phantom{Ai}  &\mbox{iff}&
(e)_0\in\VV\;\wedge\;(e)_1 \rel\phi((e)_0).
\end{eqnarray*}
 \end{definition}

 \begin{thm}\label{Theorem1a}
 $\varphi(v_1,\ldots, v_r)$ be a formula of set theory with at
 most the free variables exhibited.
 If $$\CZF+\LPO+\RDC\vdash \varphi(v_1,\ldots, v_r)$$ then
 one can explicitly construct (an index of) a partial $E$-recursive function $f$ from that  proof
  such that
  $$\BI\vdash \forall \alpha_1,\ldots,\alpha_r\in\VV\;f(\alpha_1,\ldots,\alpha_r)\rel \varphi(\alpha_1,\ldots,\alpha_r).$$

      \end{thm}
      \prf Realizability of the axioms of $\CZF+\RDC$ is just a special case of realizability over
      an $\omega\mbox{-PCA}^+$ as described in \cite[Theorem 8.5]{jucs} and is closely related to
      Aczel's \cite{aczel82} interpretation of $\CZF+\RDC$ in type theory and the realizability interpretations of $\CZF+\RDC$ presented in \cite{type94,markt,tupailo}.
      Note that to ensure realizability of $\Delta_0$ separation it is necessary that all types in $\UU$ correspond
      to sets (Corollary \ref{set}).

      We shall thus only address the realizability of $\LPO$. To avoid the niceties involved in coding functions in set theory,
      we shall demonstrate realizability of a more general type of statement which implies $\LPO$ on the basis of $\CZF$:
      $$(*)\;\;\;(\forall x\in\omega)[P(x)\,\vee\, R(x)]\to [(\exists x\in\omega)P(x)\,\vee\,(\forall x\in\omega) R(x)].$$
      To see that $(*)$ implies $\LPO$ assume that  $f\in 2^{\N}$. Then let $P(x)$ and $R(x)$ stand for $f(x)=1$ and
      $f(x)=0$, respectively.

Arguing in $\BI$, we want to show that $(*)$ is realizable.   The first step is to single out the element of $\VV$ that plays the role of the natural numbers in $\VV$. By the recursion theorem for $E$-computability define a function
$g:\N\to \N$ with index $d$ by $\{d\}^E(0)=\sup(0_{\N},\lambda x.x)$ and
    $$\{d\}^E(n+1)=\sup ((n+1)_{\N}, d\restriction n)$$ where
      $d\restriction n$ is an index of the function $g_n:\N\to \N$ with
        $g_n(k)= \{d\}^E(k)$ if $k\leq n$ and $g_n(k)=0$ otherwise.
        Finally, let $$\mathbf \omega = \sup(\natt,d).$$
        Then $\omega\in\VV$ and $\omega$ realizably plays the role of the natural numbers in $\VV$.

        Now assume that \begin{eqnarray} \label{o1} &&e\rel (\forall x\in \omega)[P(x)\,\vee\,R(x)].\end{eqnarray}
         Unraveling the definition of (\ref{o1}) we get
         $(\forall i\in \bar{\omega})\{e\}^E(i)\rel P(\tilde{\omega}i)\,\vee\,R(\tilde{\omega}i)$, whence

\begin{eqnarray} \label{o2} &&(\forall n\in \N)\{e\}^E(n)\rel P(\tilde{\omega}n)\,\vee\,R(\tilde{\omega}n).\end{eqnarray}
From (\ref{o2}) we get that for all $n\in\N$,
\begin{eqnarray} \label{o3} [(f(n))_0=0\,\wedge\,(f(n))_1\rel P(\tilde{\omega}n)]&\vee&[(f(n))_0=1\,\wedge\,(f(n))_1\rel R(\tilde{\omega}n)],\end{eqnarray}
where $f(n)=\{e\}^E(n)$. There is an index $b$ such that $\{b\}^E(n,x)= (f(n))_0$ for all $n,x$.
If there exists $n$ such that $(f(n))_0=0$ then by clause (3.2) of Definition \ref{Erec}
we get $\{\langle 3,1,b\rangle\}^E(0)=n_0+1$ where $n_0$ is the smallest number such that $(f(n_0))_0=0$.
Otherwise, by clause (3.1) of Definition \ref{Erec}, we have  $\{\langle 3,1,b\rangle\}^E(0)=0$.
We also find an index $c$ such that $\{c\}^E(k)=(n,(f(n))_1)$ if $k=n+1$ for some $n$ and
$\{c\}^E(k)=\lambda x.(f(x))_1$ if $k=0$. Let $\mathrm{sg}$ be the primitive recursive function with $\mathrm{sg}(n+1)=1$ and
$\mathrm{sg}(0)=0$. Then we have
\begin{eqnarray} (\mathrm{sg}(\{\langle 3,1,b\rangle\}^E(0)),\{c\}^E(\{\langle 3,1,b\rangle\}^E(0)))
&\rel&(\exists x\in\omega)P(x)\,\vee\,(\forall x\in\omega) R(x). \end{eqnarray}
Since there is an index $b^*$ such that $\{b^*\}^E(e)\simeq(\mathrm{sg}(\{\langle 3,1,b\rangle\}^E(0)),\{c\}^E(\{\langle 3,1,b\rangle\}^E(0)))$ this ensures the realizability of $(*)$. \qed

\paragraph{Acknowledgement:}
This material is based upon work
 supported by the EPSRC of the UK through
 Grant No. EP/G029520/1.

\end{document}